\documentclass{amsart}
\usepackage[
  left=1.0in,    
  right=1.0in,
  top=1.0in,
  bottom=1.0in,
  includefoot    
]{geometry}
\usepackage{amsthm,amsfonts,amsmath,amssymb}
\usepackage[colorlinks=true]{hyperref}
\usepackage[backrefs,msc-links,nobysame]{amsrefs}
\usepackage{cleveref}
\usepackage{tikz-cd}
\usepackage{tikz}
\usetikzlibrary{cd}
\usetikzlibrary[arrows.meta]
\usetikzlibrary{decorations.markings}
\tikzset{degil/.style={
            decoration={markings,
            mark= at position 0.5 with {
                  \node[transform shape] (tempnode) {$\backslash$};
                  }
              },
              postaction={decorate}
}
}

\DeclareMathOperator{\lann}{l.ann}
\DeclareMathOperator{\rann}{r.ann}
\DeclareMathOperator{\im}{im}

\DeclareMathOperator{\morp}{End}

\DeclareMathOperator{\cent}{Cent}

\theoremstyle{plain}
\newtheorem{theorem}{Theorem}[section] 
\newtheorem{lemma}[theorem]{Lemma}
\newtheorem{proposition}[theorem]{Proposition}
\newtheorem{corollary}[theorem]{Corollary}

\theoremstyle{definition}

\theoremstyle{remark}
\newtheorem{remark}[theorem]{Remark} 
\newtheorem{note}[theorem]{Note}

\usepackage{lineno}

\begin{document}


\pagenumbering{arabic}

\title{Central Quasi-Morphicity, Central Morphicity, and Strongly $\pi$-Regularity}

\author[Gera]{Theophilus Gera} 
\address{(Gera) Department of Mathematics, Sardar Vallabhbhai National Institute of Technology, Surat, Gujarat, India-395007}
\email{geratheophilus@gmail.com}
\thanks{The first author would like to thank the Department of Education, Government of India, for the financial assistance.}

\author[Sharma]{Amit Sharma} 
\address{(Sharma) Department of Mathematics, Sardar Vallabhbhai National Institute of Technology, Surat, Gujarat, India-395007}
\email{amitsharma@amhd.svnit.ac.in}

\date{\today}

\begin{abstract}
    This paper revisits the notions of centrally morphic and centrally quasi-morphic modules introduced by Dehghani and Sedaghatjoo. We identify sufficient conditions under which these notions become equivalent. In particular, if a module \( M \) is image-projective, generates its kernels, and its endomorphism ring has von Neumann regular center, then central morphicity, central quasi-morphicity, and right central morphicity of the endomorphism ring are equivalent. We also prove that centrally morphic rings with von Neumann regular center and ACC on principal right annihilators are strongly \(\pi\)-regular.
\end{abstract}



\maketitle

\section{Introduction}

The theory of morphic and quasi-morphic rings and modules explores how annihilator, kernel, and image relations capture the internal structure of rings and modules. Originating from the work of Nicholson and Campos~\cite{nicholson2004rings}, the notion of a \emph{morphic ring} was introduced as a dual analogue of the isomorphism theorem. Specifically, a ring $R$ is \emph{left morphic} if, for every $a\in R$, there exists an isomorphism $R/Ra\cong \lann_R(a)$ where \(\lann_R(a)\) denotes the left annihilator of \(a\). Camillo and Nicholson~\cite{camillo2007quasi} later introduced the concept of a \emph{quasi-morphic ring}, requiring only the existence of elements \(b, c \in R\) such that \(\lann_R(a) = Rb\) and \(Ra = \lann_R(c)\); the equality \(b=c\) recovers morphicity. These ideas were subsequently extended to modules by Nicholson and Campos~\cite{nicholson2005morphic} and further generalized to the quasi-morphic setting by An, Nam, and Tung~\cite{an2016quasi}. A module $M$ over a ring $R$ is said to be \emph{quasi-morphic} if, for every endomorphism $f\in \morp_R(M)$, there exist endomorphisms $g,h\in \morp_R(M)$ such that $\im f=\ker g$ and $\ker f=\im h$; if one may choose $g=h$, then $M$ is \emph{morphic}.

A more refined notion was recently proposed by Dehghani and Sedaghatjoo~\cite{dehghani2025centrally}, who introduced \emph{centrally morphic} and \emph{centrally quasi-morphic} modules by requiring the witness endomorphisms to lie in the center \(\cent(\morp_R(M))\). This additional constraint ties morphic behavior to the internal central structure of the endomorphism ring, giving rise to the following chain of implications:

\begin{center}
\begin{tikzcd}[row sep=small, column sep=small]
    & & \text{morphic} \arrow[dr,Rightarrow] & \\
    \text{strongly endoregular} \arrow[r,Rightarrow] & \text{centrally morphic} \arrow[ur,Rightarrow] \arrow[dr,Rightarrow] & & \text{quasi-morphic} \\
    & & \text{centrally quasi-morphic} \arrow[ur,Rightarrow] &
\end{tikzcd}
\end{center}

Dehghani and Sedaghatjoo established in \cite{dehghani2025centrally}*{Proposition 2.2} that for image-projective modules, central morphicity and central quasi-morphicity coincide. They also showed that projective modules, being image-projective, satisfy this equivalence (Corollary 2.3). Their work left open the general question of whether the equivalence holds for arbitrary modules (Question 2.14).

The present paper identifies additional sufficient conditions under which these notions become equivalent. We show that if $M$ is image-projective, generates its kernels, and $Z(S)$ is von Neumann regular, then the following are equivalent: $M$ is centrally morphic; $M$ is centrally quasi-morphic; and $S=\morp_R(M)$ is right centrally morphic (see Proposition \ref{prop:corrected}). This result extends \cite{dehghani2025centrally}*{Proposition 2.2} by adding the von Neumann regularity condition on the center and proving the equivalence with right central morphicity of $S$. For projective modules that generate their kernels, the equivalence reduces to: $M$ is centrally morphic if and only if $M$ is centrally quasi-morphic and $S$ is right centrally morphic, provided $Z(S)$ is von Neumann regular (see Corollary \ref{cor:proj-central}).

From a ring-theoretic viewpoint, these results further connect central morphicity with strong $\pi$-regularity, providing a different perspective on the classical proposition of Lee and Zhou~\cite{lee2009regularity}*{Proposition 23}. We show that a centrally morphic ring whose center is von Neumann regular and which satisfies ACC on principal right annihilators must itself be strongly $\pi$-regular (see Theorem~\ref{thm:cqm-acc}). In summary, this paper provides sufficient conditions for the equivalence of central morphicity and central quasi-morphicity, extending the framework of Dehghani and Sedaghatjoo and clarifying the role of centrality and regularity in morphic module theory. Section~\ref{sec:results} presents the main results, including criteria for the equivalence between central morphicity and central quasi-morphicity, along with structural consequences for projective modules and their endomorphism rings.

Throughout, all rings are associative with identity, and all modules are unital right modules unless stated otherwise. For a module $M$, we denote its endomorphism ring by $\morp_R(M)$ (or $S$), acting on the left. The left and right annihilators in $R$ are denoted by $\lann_R$ and $\rann_R$, respectively. Standard references include \citelist{\cite{lam2001firstcourse} \cite{lam1999lectures}}.

\section{Main Results} \label{sec:results}

A module $M$ is called \emph{image-projective} if for any $f, g \in S = \morp_R(M)$, the inclusion $\im(f) \subseteq \im(g)$ implies $f \in gS$. In other words, containment among images of endomorphisms is reflected algebraically within the ring $\morp_R(M)$. Image-projective modules were introduced by Nicholson and Campos~\cite{nicholson2005morphic} in their study of morphic module structures.

A module $M$ \emph{generates its kernels} if for every endomorphism $f \in S = \morp_R(M)$, the submodule $\ker f$ is generated by images of endomorphisms of $M$ contained in it; that is,
\[
\ker f = \sum \{ \im(\gamma) \mid \gamma \in S, \ \im(\gamma) \subseteq \ker f \}.
\]
Equivalently, for each $x \in \ker f$, there exist finitely many $\gamma_1, \dots, \gamma_n \in S$ with $\im(\gamma_i) \subseteq \ker f$ such that $x \in \sum_i \im(\gamma_i)$.

The following lemma records two useful properties of image-projective modules that will be used repeatedly.

\begin{lemma} \label{lm:imageprojective}
    Let $M$ be an $R$-module and $S = \morp_R(M)$.
    \begin{enumerate}
        \item If $\ker f = \im g$ for some $g \in S$, then $M$ generates $\ker f$.
        \item Suppose $M$ generates its kernels and $e \in S$ is an idempotent satisfying $\rann_S(f) = eS$. Then $\ker f = eM$.
    \end{enumerate}
\end{lemma}

\begin{proof}
    \begin{enumerate}
        \item Immediate, since $\ker f = \im g$ is by definition generated by the image of the map $g : M \to M$.
        \item Because $fe = 0$, we have $eM \subseteq \ker f$. For the converse, let $x \in \ker f$. Since $M$ generates its kernels, there exist finitely many $\gamma_i \in S$ with $\im \gamma_i \subseteq \ker f$ and $x \in \sum_i \im \gamma_i$. Each such $\gamma_i$ satisfies $f \gamma_i = 0$, hence $\gamma_i \in \rann_S(f) = eS$; say $\gamma_i = e s_i$ for some $s_i \in S$. Then $\im \gamma_i \subseteq eM$, so $x \in eM$. Thus $\ker f = eM$.
    \end{enumerate}
\end{proof}

We now turn to the interaction between module-theoretic and ring-theoretic structures. A ring $R$ is said to be \emph{right centrally quasi-morphic} if for each $a \in R$, there exist central elements $b, c \in Z(R)$ such that $\rann_R(a) = bR$ and $aR = \rann_R(c)$. If the same holds on both sides, $R$ is called \emph{centrally quasi-morphic}. If $b=c$, then $R$ is called \emph{centrally morphic}. The following result extends \cite{dehghani2025centrally}*{Proposition 2.2} by adding the hypothesis that $Z(S)$ is von Neumann regular and proving the equivalence with right central morphicity of $S$.

\begin{proposition} \label{prop:corrected}
Let $M$ be a right $R$-module and $S = \morp_R(M)$. Suppose that
\begin{itemize}
    \item $M$ is image‑projective,
    \item $M$ generates its kernels,
    \item $Z(S)$ is von Neumann regular.
\end{itemize}
Then the following are equivalent:
\begin{enumerate}
    \item $M$ is centrally morphic;
    \item $M$ is centrally quasi‑morphic;
    \item $S$ is right centrally morphic.
\end{enumerate}
\end{proposition}

\begin{proof}
{(1) $\Rightarrow$ (2)}. Immediate.

{(2) $\Rightarrow$ (3)}. Assume $M$ is centrally quasi‑morphic. For each $f \in S$ there exist central elements $g, h \in Z(S)$ such that $\ker f = \im g$ and $\im f = \ker h$. Because $M$ is image‑projective, $\im \varphi \subseteq \ker f = \im g$ implies $\varphi \in gS$. Hence
\begin{align}
    \rann_S(f) = gS \label{prop2.2:1}
\end{align}

We now prove $fS = \rann_S(h)$.

($\subseteq$). Since $\im f = \ker h$, for every $x \in M$ we have $f(x) \in \ker h$, so $h(f(x)) = 0$. Thus $hf = 0$. For any $s \in S$, $h(fs) = (hf)s = 0$; therefore $fS \subseteq \rann_S(h)$.

($\supseteq$). Let $\psi \in \rann_S(h)$. Then $h\psi = 0$, so for all $x \in M$, $h(\psi(x)) = 0$; hence $\psi(x) \in \ker h = \im f$. Thus $\im\psi \subseteq \im f$. By image‑projectivity, $\psi \in fS$.

Consequently,
\begin{align}
    fS = \rann_S(h) \label{prop2.2:2}
\end{align}
Thus $S$ is right centrally quasi‑morphic.

Because $Z(S)$ is von Neumann regular, for any central element $z$ there exists $x \in Z(S)$ with $z = zxz$. Set $e = zx$; then $e^2 = e$, $e \in Z(S)$, and $zS = eS$. Apply this to $g$ and $h$ to obtain central idempotents $e, e' \in Z(S)$ such that
\begin{align}
    gS = eS,\qquad hS = e'S \label{prop2.2:3}
\end{align}

From \eqref{prop2.2:1} and \eqref{prop2.2:3} we have $\rann_S(f) = eS$. From \eqref{prop2.2:2} we have $fS = \rann_S(h)$; we now compute $\rann_S(h)$.

\emph{Claim}. $\rann_S(h) = (1-e')S$.

\begin{proof}
Since $hS = e'S$, there exist $u,v \in S$ with $h = e'u$ and $e' = hv$. Because $h \in Z(S)$, we have $hv = vh$; therefore $e' = hv = vh$. For any $x \in S$:
\begin{itemize}
    \item If $e'x = 0$, then using centrality of $e'$ ($e'u = ue'$): $hx = e'ux = ue'x = 0$.
    \item If $hx = 0$, then $e'x = vhx = 0$.
\end{itemize}
Thus $hx = 0 \iff e'x = 0$. Hence $\rann_S(h) = \rann_S(e')$. Since $e'$ is a central idempotent, $\rann_S(e') = (1-e')S$.
\end{proof}

Therefore, from \eqref{prop2.2:2} we obtain
\begin{align}
    fS = (1-e')S \label{prop2.2:4}
\end{align}
We now show that $e = e'$. From $\rann_S(f) = eS$ we have $fe = 0$; central idempotents commute, so $ef = 0$. Using \eqref{prop2.2:4}, write $f = (1-e')u$ and $1-e' = fv$ for some $u,v \in S$. Then
\begin{align}
    e(1-e') = efv = 0 \;\Longrightarrow\; e = ee' \label{prop2.2:5}
\end{align}

Also, $e'f = e'(1-e')u = 0$; centrality of $e'$ gives $fe' = 0$. Hence $e' \in \rann_S(f) = eS$, so $e' = es$ for some $s \in S$. Multiplying on the left by $(1-e)$ yields $(1-e)e' = 0$, i.e. $e' = ee'$. Together with \eqref{prop2.2:5} we obtain $e = e'$.

Consequently, $fS = (1-e)S = \rann_S(e)$. Thus $S$ is right centrally morphic.

{(3) $\Rightarrow$ (1)}. Assume $S$ is right centrally morphic. For any $f \in S$ there exists a central idempotent $e \in Z(S)$ such that $\rann_S(f) = eS$ and $fS = \rann_S(e)$. Because $M$ generates its kernels, Lemma~\ref{lm:imageprojective}(2) gives $\ker f = eM$.

Since $e$ is a central idempotent, $\rann_S(e) = (1-e)S$; hence $fS = (1-e)S$. Thus there exist $u,v \in S$ with $f = (1-e)u$, and $1-e = fv$. Now $(1-e)M = fv(M) \subseteq fM$, so $(1-e)M \subseteq \im f$. Conversely, for any $x = f(y) \in \im f$, $x = f(y) = (1-e)u(y) \in (1-e)M$. Therefore $\im f = (1-e)M$.

Thus for every $f \in S$ there exists a central idempotent $e \in Z(S)$ such that $\ker f = eM$ and $\im f = (1-e)M$, which means $M$ is centrally morphic.
\end{proof}

The following result refines the correspondence between morphic and quasi-morphic structures in the semisimple artinian setting. It extends \cite{gera2025modules-2}*{Corollary~3.21} from finitely generated modules over semisimple artinian rings to arbitrary semisimple artinian modules, while making explicit the role of the endomorphism ring in mediating the equivalence between module-theoretic and ring-theoretic centrality conditions.

\begin{corollary}\label{cor:semisimple}
Let $M$ be a semisimple artinian right $R$-module and set $S=\morp_R(M)$. 
Then the following are equivalent:
\begin{enumerate}
    \item $M$ is centrally quasi-morphic;
    \item $M$ is centrally morphic;
    \item $S$ is right centrally morphic.
\end{enumerate}
Moreover, $Z(S)$ is a finite direct product of fields and hence von Neumann regular (indeed strongly regular).
\end{corollary}

\begin{proof}
Since $M$ is semisimple, every submodule is a direct summand. It follows that $M$ is image-projective and generates its kernels. Moreover, $S=\morp_R(M)$ is semisimple artinian, hence  $Z(S)$ is a finite product of fields and therefore von Neumann regular. Thus the hypotheses of Proposition~\ref{prop:corrected} are satisfied, yielding the stated equivalence.

Finally, the structure theorem $S\cong \prod_i M_{\alpha_i}(D_i)$ implies 
$Z(S)\cong \prod_i Z(D_i)$, a finite product of fields, which is von Neumann regular.
\end{proof}

The following corollary provides sufficient conditions under which the equivalence holds for projective modules.

\begin{corollary}\label{cor:proj-central}
Let $M$ be a projective right $R$-module that generates its kernels. Set $S=\morp_R(M)$ and assume that $Z(S)$ is von Neumann regular. Then the following are equivalent:
\begin{enumerate}
    \item $M$ is centrally morphic;
    \item $M$ is centrally quasi-morphic;
    \item $S$ is right centrally morphic.
\end{enumerate}
\end{corollary}

\begin{proof}
If $M$ is projective then it is image‑projective. Together with the hypothesis, Proposition~\ref{prop:corrected} applies directly, yielding the stated equivalence.
\end{proof}

\begin{note}
Whether every centrally quasi-morphic module is centrally morphic remains an open problem in general (see Question 2.14 of \cite{dehghani2025centrally}). Our results provide sufficient conditions under which the equivalence holds, notably for modules satisfying the hypotheses of Proposition~\ref{prop:corrected} and its corollaries. We also note that centrally morphic modules need not be Rickart (see, for example, \( R = \Bbbk[x]/(x^2) \) as a module over itself), but this does not settle the general question.
\end{note}

Motivated by Proposition~23 of Lee and Zhou~\cite{lee2009regularity}, which characterizes strongly regular rings via semiprimeness and central morphicity, we investigate what additional hypotheses beyond central morphicity yield strong $\pi$-regularity. While semiprimeness yields strong regularity (and hence strong $\pi$-regularity) in the result of Lee and Zhou, our approach instead uses stabilization of annihilator decompositions associated to powers of elements. We show that this stabilization follows from the ascending chain condition on principal right annihilators together with the assumption that the center is von Neumann regular. Under these conditions, every centrally morphic ring is strongly $\pi$-regular.

\begin{theorem}\label{thm:cqm-acc}
Let $R$ be a centrally morphic ring. Suppose that:
\begin{enumerate}
    \item $Z(R)$ is von Neumann regular;
    \item $R$ satisfies the ascending chain condition on principal right annihilators.
\end{enumerate}
Then $R$ is strongly $\pi$-regular.
\end{theorem}

\begin{proof}
    Let $a\in R$. Because $R$ is centrally morphic, for every positive integer $k$, there exist central elements $b_k,c_k \in Z(R)$ such that $\rann_R(a^k)=b_k R$, $\rann_R(b_k)=a^k R$, $\lann_R(a^k)=Rc_k$ and $\lann_R(c_k)=Ra^k$.

    Since $Z(R)$ is von Neumann regular, for each $k$ we can choose $x_k \in Z(R)$ with $b_k =b_k x_k b_k$; then $e_k \in Z(R), e_k^2=e_k$ and $b_k R=e_k R$. Consequently,
    \begin{align}
        \rann_R(a^k)=e_k R, \quad a^k R=\rann_R(e_k)=(1-e_k)R \label{thm2.6:1}
    \end{align}

    Similarly, using the regularity of $Z(R)$ on $c_k$, we obtain central idempotents $f_k \in Z(R)$ such that 
    \begin{align}
        \lann_R(a^k)=Rf_k, \quad Ra^k=\lann_R(f_k)=R(1-f_k) \label{thm2.6:2}
    \end{align}

    Now we analyze the right-hand side. From $a^{k+1} R\subseteq a^k R$ and \eqref{thm2.6:1} we have $(1-e_{k+1})R \subseteq (1-e_k) R$. Multiplying on the right by $e_k$ gives $(1-e_{k+1})e_k R=0$, hence $(1-e_{k+1})e_k=0$. Thus, $e_k=e_{k+1} e_k$, and because the idempotents are central, $e_k=e_k e_{k+1}$. Therefore, $e_k \in e_{k+1} R$ and consequently $e_k R\subseteq e_{k+1} R$. Hence $\{e_k R \}$ is an ascending chain of principal right annihilators (each $e_k R$ equals $\rann_R(1-e_k)$). By hypothesis this chain stabilizes: there exists $N$ such that $e_k R=e_N R$ for all $k\geq N$.

    Since each $e_k$ is a central idempotent, equality of the principal right ideals forces equality of the idempotents themselves; thus $e_k=e_N$ for all $k\geq N$. From \eqref{thm2.6:1}, we obtain $a^k R=(1-e_N)R$ for all $k\geq N$, and in particular 
    \begin{align}
        a^N R=a^{N+1} R \label{thm2.6:4}
    \end{align}

    Now treat the left-hand side. From $Ra^{k+1} \subseteq Ra^k$ and \eqref{thm2.6:2}, we have $R(1-f_{k+1})\subseteq R(1-f_k)$. For any $x\in R$, write $x(1-f_{k+1})=y(1-f_k)$ with some $y\in R$. Multiplying on the right by $f_k$ yields $x(1-f_{k+1})f_k=y(1-f_k)f_k=0$. Since this holds for every $x$, we get $(1-f_{k+1})f_k=0$, hence $f_k=f_{k+1} f_k$. Centrality of the idempotents gives $f_k =f_k f_{k+1}$, so $f_k \in f_{k+1} R$ and therefore $f_k R\subseteq f_{k+1} R$. Thus $\{f_k R\}$ is again an ascending chain of principal right annihilators (because $f_k R=\rann_R(1-f_k)$). By the ACC there exists $M$ such that $f_k R=f_M R$ for all $k\geq M$. As before, equality of the right ideals implies $f_k =f_M$ for all $k\geq M$. From \eqref{thm2.6:2} we obtain $Ra^k=R(1-f_M)$ for all $k\geq M$, and in particular

    \begin{align}
        Ra^M=Ra^{M+1} \label{thm2.6:6}
    \end{align}

    Let $t=\max \{N,M \}$. Then from \eqref{thm2.6:4} and \eqref{thm2.6:6}, we have $a^t R=a^{t+1} R$, and $Ra^t=Ra^{t+1}$. Hence there exists elements $r,u \in R$ such that $a^t=a^{t+1} r$ and $a^t=ua^{t+1}$. Since $a$ was arbitrary, $R$ is strongly $\pi$-regular.
\end{proof}

\begin{remark}
\cite{lee2009regularity}*{Proposition~23} shows that every semiprime centrally morphic ring is strongly regular, and hence nonsingular. Consequently, the open question ``Are semiprime $\pi$-regular rings nonsingular?'' \cite{hongying2026regularity}*{Question, p. 6} has a positive answer within the subclass of semiprime centrally morphic rings. 
\end{remark}

\bibliography{ref}
\end{document}